\documentclass[p1]{elsarticle}
\usepackage[utf8]{inputenc}
\usepackage{amsfonts}
\usepackage{amssymb}
\usepackage{amsmath}
\usepackage{hyperref}
\hypersetup{colorlinks = true, citecolor = blue, linkcolor = blue}

\title{On the properties of falling and rising factorial transforms}
\author{Parham Zarghami\corref{cor1}\fnref{fn1}} 
\cortext[cor1]{Corresponding author \ead{parham.zarghami@ut.ac.ir}}
\fntext[fn1]{Declarations of interest: None}

\affiliation{organization={University of Tehran}, addressline={School of Electrical and Computer Engineering, University College of Engineering, North Kargar street}, city={Tehran}, postcode={1439957131}, country={Iran}}

\date{December 2023}

\begin{document}

\begin{abstract}
    In this paper, we will extend the falling and rising factorial transforms \cite{ref. 1} which in this case every arbitrary function can be applied. Then, the properties of these transforms will be investigated and some corollaries will be shown.
    
    These transforms have interesting properties that led to new series expansions, representation of functions in terms of operators, connection between different series and operators, etc. By converting the power series into the Newton series, useful identities can be derived.
\end{abstract}


\maketitle
\section{Introduction}
\label{section 1}

The falling and rising factorial transforms introduced in \cite{ref. 1} have been utilized to convert the power series into factorial ones. This conversion causes the discovery of hidden aspects of some intrinsically combinatorial functions (e.g., gamma function) that cannot be observed by known series.

Based on the definition of the falling and rising factorial transforms and their inverses \cite{ref. 1}, we can represent the following generalized formulas:
\begin{equation}
	\label{eq. 1}
	FFT (\sum_{n=0}^{\infty} a_n x^n) = \sum_{n=0}^{\infty} a_n (x)_n
\end{equation}
\begin{equation}
	\label{eq. 2}
	FFT^{-1} (\sum_{n=0}^{\infty} a_n (x)_n) = \sum_{n=0}^{\infty} a_n x^n
\end{equation}
\begin{equation}
	\label{eq. 3}
	RFT (\sum_{n=0}^{\infty} a_n x^n) = \sum_{n=0}^{\infty} a_n x^{\overline{n}}
\end{equation}
\begin{equation}
	\label{eq. 4}
	RFT^{-1} (\sum_{n=0}^{\infty} a_n x^{\overline{n}}) = \sum_{n=0}^{\infty} a_n x^n
\end{equation}
Which $a_n=\partial_x^n f(0)/n!$ or $a_n=\Delta_x^n f(0)/n!$ can be considered.$x^{\overline{n}}$ and $(x)_n$ are rising and falling factorials, respectively. Also, $\partial_x$ and $\Delta_x$ are derivative and difference symbols on $x$. (\ref{eq. 1}), (\ref{eq. 2}), (\ref{eq. 3}), and (\ref{eq. 4}) demonstrate the definition of FFT (i.e., falling factorial transform) and RFT (i.e., rising factorial transform) and are not practical. Other practical formulations suggested by \cite{ref. 1} that we can generalize them as follows:
\begin{equation}
	\label{eq. 5}
	FFT (f(x)) = \partial_t^x (f(t) e^t ) \Big\rvert_{t=0}
\end{equation}
\begin{equation}
	\label{eq. 6}
	FFT^{-1} (f(x)) = e^{-x} \sum_{n=0}^{\infty} \frac{f(n)}{n!} x^n
\end{equation}
\begin{equation}
	\label{eq. 7}
	RFT (f(x)) = \frac{1}{\Gamma (x)} \int_{0}^{\infty} f(t) t^{x-1} e^{-t} dt
\end{equation}
\begin{equation}
	\label{eq. 8}
	RFT (f(x)) = \frac{1}{\Gamma (x)}  \mathcal{M} \Bigl(f(x)e^{-x} \Bigr)
\end{equation}
\begin{equation}
	\label{eq. 9}
	RFT^{-1} (f(x)) = e^x \sum_{n=0}^{\infty} \frac{(-1)^n f(-n)}{n!} x^n
\end{equation}
Bear in mind that variables before and after transformation must be determined precisely (e.g., $FFT_x (f(x))(s)$). But in some equations like the above, it is not necessary to write the right side of the equation in terms of other variables.

Table \ref{tab:Table 1} shows the properties of Mellin transform \cite{ref. 6}.
\begin{table}[h!]
	\caption{Properties of Mellin transform}
	\label{tab:Table 1}
	\centering
	\begin{tabular}{|c|c|}
		\hline
		\textbf{Function} & \textbf{Mellin transform} \\
		\hline
		$x^a f(x)$ & $\mathcal{M} \Bigl(f(x) \Bigr)(x+a)$ \\
		\hline
		$f(ax)$ & $a^{-x} \mathcal{M}(f(x))$ \\
		\hline
		$f(x^a)$ & $\lvert a \rvert^{-1} \mathcal{M} \Bigl(f(x) \Bigr)(x/a)$ \\
		\hline
		$f(x)(log(x))^n$ & $\partial_x^n \mathcal{M}(f(x))$ \\
		\hline
		$\partial_x^n f(x)$ & $(-1)^n  \frac{\Gamma (x)}{\Gamma (x-n)} \mathcal{M} \Bigl(f(x) \Bigr)(x-n)$ \\
		\hline
		$\int_{0}^{x} f(y)dy$ & $\frac{-1}{x} \mathcal{M} \Bigl(f(x) \Bigr)(x+1)$ \\
		\hline
		$\int_{x}^{\infty} f(y)dy$ & $\frac{1}{x} \mathcal{M} \Bigl(f(x) \Bigr)(x+1)$ \\
		\hline
		$(x\partial_x)^n f(x)$ & $(-1)^n x^n \mathcal{M} \Bigl(f(x) \Bigr)$ \\
		\hline
		$(\partial_x x)^n f(x)$ & $(-1)^n (x-1)^n \mathcal{M}(f(x))$ \\
		\hline
		$f(x)g(x)$ & $\frac{1}{2 \pi i} \int_{c-i\infty}^{c+i\infty} \mathcal{M} \Bigl(f(x) \Bigr)(r) \mathcal{M} \Bigl(g(x) \Bigr)(x-r)dr$ \\
		\hline
		$\int_{0}^{\infty} f(x/y)g(y)dy/y$ & $\mathcal{M}(f(x))\mathcal{M}(g(x))$ \\
		\hline
	\end{tabular}
\end{table}

In section \ref{section 2}, we will investigate the properties of falling and rising factorial transforms. Then in section \ref{section 3}, some examples will be shown. Finally, future works and conclusions will be discussed.

\section{Properties}
\label{section 2}

Most of the properties of FFT and RFT are the same with slight differences in formulas. So concisely, we will investigate these similar properties only for the falling factorial transform. This similarity of the mentioned transforms originated from the following formula:
\begin{equation}
	\label{eq. 10}
	RFT_t \Bigl(f(t) \Bigr)(x) = FFT_t \Bigl(f(-t) \Bigr)(-x)
\end{equation}
To prove the above equation, we use (\ref{eq. 5}) and \cite{ref. 1}:
\begin{equation}
	\label{eq. 11}
	\begin{split}
		FFT_t \Bigl(f(-t) \Bigr)(-x) & = \partial_t^{-x} (f(-t) e^t) \Big \rvert_{t=0} = \frac{1}{\Gamma (x)} \int_{-\infty}^{0} f(-t) e^t (-t)^{x-1} dt \\
		& = \frac{1}{\Gamma (x)} \int_{0}^{\infty} f(t) t^{x-1} e^{-t} dt = RFT_t \Bigl(f(t) \Bigr)(x)
	\end{split} 
\end{equation}
This occurs due to $\sum_{n \geq 0} a_n x^{\overline{n}} = \sum_{n \geq 0} a_n (-1)^n (-x)_n$. (\ref{eq. 10}) is also valid for inverse transforms owing to (\ref{eq. 6}) and (\ref{eq. 9}). Hence by variable changing or using the same procedure, the following formulas can be derived for RFT.

Falling factorial transform has the following properties:

\subsection{Linearity}

By (\ref{eq. 5}) and (\ref{eq. 6}), we can write:

\begin{equation}
	\label{eq. 12}
	\begin{split}
		& FFT(\alpha f(x) + \beta g(x)) = \partial_t^x ((\alpha f(t) + \beta g(t))e^t) \Big\rvert_{t=0} \\
		& = \alpha \partial_t^x (f(t) e^t) \Big\rvert_{t=0} + \beta \partial_t^x (g(t) e^t) \Big\rvert_{t=0} = \alpha FFT(f(x)) + \beta FFT(g(x))
	\end{split} 
\end{equation}
\begin{equation}
	\label{eq. 13}
	\begin{split}
		& FFT^{-1} (\alpha f(x) + \beta g(x)) = e^{-x} \sum_{n=0}^{\infty} \frac{\alpha f(n) + \beta g(n)}{n!} x^n = \alpha e^{-x} \sum_{n=0}^{\infty} \frac{f(n)}{n!} x^n \\
		& + \beta e^{-x} \sum_{n=0}^{\infty} \frac{g(n)}{n!} x^n = \alpha FFT^{-1} (f(x)) + \beta FFT^{-1} (g(x))
	\end{split}
\end{equation}

\subsection{Derivatives and differences}
\label{subsection 2.2}

We can apply $\partial_x^k (x^n) = (n)_k x^{n-k}$ and $\Delta_x^k ((x)_n) = (n)_k (x)_{n-k}$ to (\ref{eq. 1}) and (\ref{eq. 2}) respectively:
\begin{equation}
	\label{eq. 14}
	\begin{split}
		FFT(\partial_x^k (x^n)) & = FFT((n)_k x^{n-k}) = (n)_k (x)_{n-k} = \Delta_x^k ((x)_n) \\
		& = \Delta_x^k FFT(x^n)
	\end{split}
\end{equation}
\begin{equation}
	\label{eq. 15}
    \begin{split}
    	FFT^{-1} (\Delta_x^k ((x)_n)) & = FFT^{-1} ((n)_k (x)_{n-k}) = (n)_k x^{n-k} = \partial_x^k (x^n) \\
    	& = \partial_x^k FFT^{-1} ((x)_n)
    \end{split}
\end{equation}
Generally:
\begin{equation}
	\label{eq. 16}
    FFT (\partial_x^k f(x)) = \Delta_x^k FFT (f(x))
\end{equation}
\begin{equation}
    \label{eq. 17}
    FFT^{-1} (\Delta_x^k f(x)) = \partial_x^k FFT^{-1} (f(x))
\end{equation}
From $\partial_x = log(1+ \Delta_x)$, we can obtain:
\begin{equation}
    \label{eq. 18}
    \begin{split}
        FFT^{-1} (\partial_x^k f(x)) & = FFT^{-1} ((log(1 + \Delta_x))^k f(x)) \\
        & = (log(1 + \partial_x))^k FFT^{-1} (f(x))
    \end{split}
\end{equation}
This procedure can be repeated for FFT of differences as follows:
\begin{equation}
	\label{eq. 19}
    \begin{split}
    	FFT(\Delta_x^k f(x)) & = FFT((e^{\partial_x}-1)^k f(x)) \\
    	& = (e^{\Delta_x}-1)^k FFT(f(x))
    \end{split}
\end{equation}

\subsection{Integration and summation}

(\ref{eq. 16}), (\ref{eq. 17}), (\ref{eq. 18}), and (\ref{eq. 19}) are hold for negative values of $k$. So, we can get:
\begin{equation}
	\label{eq. 20}
	FFT (\partial_x^{-k} f(x)) = \Delta_x^{-k} FFT (f(x)) 
\end{equation}
\begin{equation}
	\label{eq. 21}
    FFT^{-1} (\Delta_x^{-k} f(x)) = \partial_x^{-k} FFT^{-1} (f(x))
\end{equation}
\begin{equation}
	\label{eq. 22}
    \begin{split}
    	FFT^{-1} (\partial_x^{-k} f(x)) & = FFT^{-1} ((log(1 + \Delta_x))^{-k} f(x)) \\
    	& = (log(1 + \partial_x))^{-k} FFT^{-1} (f(x))
    \end{split}
\end{equation}
\begin{equation}
	\label{eq. 23}
    \begin{split}
    	FFT(\Delta_x^{-k} f(x)) & = FFT((e^{\partial_x}-1)^{-k} f(x)) \\
    	& = (e^{\Delta_x}-1)^{-k} FFT(f(x))
    \end{split}
\end{equation}
One can prove these equations by the procedure discussed in section \ref{subsection 2.2}. 

In the case of definite integral over $(0,\infty)$ on inverse falling factorial transform, we can obtain a useful identity by (\ref{eq. 6}):
\begin{equation}
\label{eq. 24}
    \int_{0}^{\infty} FFT^{-1} (f(t))dt = \int_{0}^{\infty} e^{-t} \sum_{n=0}^{\infty} \frac{f(n)}{n!} t^n dt = \sum_{n=0}^{\infty} \frac{f(n)}{n!} n! = \sum_{n=0}^{\infty} f(n) 
\end{equation}

\subsection{New operators}
\label{subsection 2.4}

By $(x)_n = \sum_{k=0}^{n} {n\brack k} (-1)^{n-k} x^k$ and (\ref{eq. 1}):
\begin{equation}
	\label{eq. 25}
	\begin{split}
		FFT(f(x)) & = \sum_{n=0}^{\infty} \frac{\partial_x^n f(0)}{n!} (x)_n = \sum_{n=0}^{\infty} \sum_{k=0}^{n} {n\brack k} (-1)^{n-k} \frac{\partial_x^n f(0)}{n!} x^k \\
		& = \sum_{k=0}^{\infty} \sum_{n=k}^{\infty} {n\brack k} (-1)^{n-k} \frac{\partial_x^n f(0)}{n!} x^k = \sum_{k=0}^{\infty} \frac{(log(1+\partial_x))^k f(0)}{k!} x^k
	\end{split}
\end{equation}
Also, by $x^n = \sum_{k=0}^{n} {n\brace k} (x)_k$ and (\ref{eq. 2}):
\begin{equation}
	\label{eq. 26}
	\begin{split}
		FFT^{-1} (f(x)) & = \sum_{n=0}^{\infty} \frac{\Delta_x^n f(0)}{n!} x^n = \sum_{n=0}^{\infty} \sum_{k=0}^{n} {n\brace k} \frac{\Delta_x^n f(0)}{n!} (x)_k \\
		& = \sum_{k=0}^{\infty} \sum_{n=k}^{\infty} {n\brace k} \frac{\Delta_x^n f(0)}{n!} (x)_k = \sum_{k=0}^{\infty} \frac{(e^{\Delta_x}-1)^k f(0)}{k!} (x)_k
	\end{split}
\end{equation}
Operators $(log(1+ \partial_x))^k$ and $(e^{\Delta_x}-1)^k$ were emerged in previous sections. We can consider $Z_n (x) = FFT((x)_n) = \sum_{k=0}^{n} {n\brack k} (-1)^{n-k} (x)_k$. Thus, with (\ref{eq. 18}), (\ref{eq. 19}), and $FFT^{-1} (x^n) = T_n (x)$ \cite{ref. 1}:
\begin{equation}
	\label{eq. 27}
	(log(1+\partial_x))^k T_n (x) = (n)_k T_{n-k} (x)
\end{equation}
\begin{equation}
	\label{eq. 28}
	(e^{\Delta_x}-1)^k Z_n (x)= (n)_k Z_{n-k} (x)
\end{equation}
$T_n (x)$ is the $n$th Touchard polynomial. Two new series expansions can be discovered by applying the inverse FFT and FFT on (\ref{eq. 25}) and (\ref{eq. 26}) respectively:
\begin{equation}
	\label{eq. 29}
	f(x) = \sum_{k=0}^{\infty} \frac{(log(1+\partial_x))^k f(0)}{k!} T_k (x)
\end{equation}
\begin{equation}
	\label{eq. 30}
	f(x) = \sum_{k=0}^{\infty} \frac{(e^{\Delta_x}-1)^k f(0)}{k!} Z_k (x)
\end{equation}
In general, one can rewrite these formulas as follows:
\begin{equation}
	\label{eq. 31}
	f(x) = \sum_{k=0}^{\infty} \frac{(log(1+\partial_x))^k f(x_0)}{k!} T_k (x - x_0)
\end{equation}
\begin{equation}
	\label{eq. 32}
	f(x) = \sum_{k=0}^{\infty} \frac{(e^{\Delta_x}-1)^k f(x_0)}{k!} Z_k (x - x_0)
\end{equation}

\subsection{Shifting}

We know that $f(x+a) = E^a f(x)=e^{a \partial_x} f(x)$. So:
\begin{equation}
	\label{eq. 33}
	\begin{split}
		FFT(f(x+a)) & = FFT(e^{a \partial_x} f(x)) = e^{a \Delta_x} FFT(f(x)) \\
		& = \sum_{n=0}^{\infty} \frac{a^n}{n!} \Delta_x^n FFT(f(x))
	\end{split}
\end{equation}
\begin{equation}
	\label{eq. 34}
	\begin{split}
		FFT^{-1} (f(x+a)) & = FFT^{-1} (e^{a \partial_x} f(x)) = e^{a log(1+\partial_x)} FFT^{-1} (f(x)) \\
		& = (1+\partial_x)^a FFT^{-1} (f(x)) = \sum_{n=0}^{\infty} \binom{a}{n} \partial_x^n FFT^{-1} (f(x))
	\end{split}
\end{equation}
Now we can take $FFT_a$ and $FFT_a^{-1}$ on (\ref{eq. 33}) and (\ref{eq. 34}):
\begin{equation}
	\label{eq. 35}
	\begin{split}
		& FFT_a (FFT_x (f(x+a))) = \sum_{n=0}^{\infty} \frac{(a)_n}{n!} \Delta_x^n FFT_x (f(x)) \\
		& = (1+ \Delta_x)^a FFT_x (f(x))= E^a FFT_x (f(x)) = FFT_t (f(t))(x+a)
	\end{split}
\end{equation}
\begin{equation}
	\label{eq. 36}
	\begin{split}
		FFT_a^{-1} (FFT_x^{-1} (f(x+a))) & = \sum_{n=0}^{\infty} \frac{a^n}{n!} \partial_x^n FFT_x^{-1} (f(x)) \\
		& = e^{a \partial_x} FFT_x^{-1} (f(x)) = FFT_t^{-1} \Bigl(f(t) \Bigr)(x+a)
	\end{split}
\end{equation}
In (\ref{eq. 35}) and (\ref{eq. 36}), inner and outer shifting coincided.

\subsection{Orthogonality}

Another way to investigate shifting property is to use the definition of falling factorial transform and its inverse:
\begin{equation}
	\label{eq. 37}
	\begin{split}
		FFT((x+a)^n) & = FFT(\sum_{k=0}^{n} \binom{n}{k} x^k a^{n-k}) = \sum_{k=0}^{n} \binom{n}{k} (x)_k a^{n-k} \\
		& = n! L_n^{(x-n)} (-a) = (-a)^n c_n (x,-a)
	\end{split}
\end{equation}
\begin{equation}
	\label{eq. 38}
	\begin{split}
		FFT^{-1} ((x+a)_n) & = FFT^{-1} (\sum_{k=0}^{n} \binom{n}{k} (x)_k (a)_{n-k}) = \sum_{k=0}^{n} \binom{n}{k} x^k (a)_{n-k} \\
		& = n! L_n^{(a-n)} (-x) = (-x)^n c_n (a,-x)
	\end{split}
\end{equation}
Which $L_n^{(\alpha)} (x)$ and $c_n (x,a)$ are generalized Laguerre polynomials and Charlier polynomials, respectively \cite{ref. 3}. (\ref{eq. 37}) and (\ref{eq. 38}) are obvious concerning the definition of mentioned polynomials. Charlier polynomials satisfy the below orthogonality equation \cite{ref. 3}:
\begin{equation}
	\label{eq. 39}
	\sum_{k=0}^{\infty} \frac{a^k}{k!} c_n (k,a) c_m (k,a) = \frac{e^a n!}{a^n} \delta_{n,m}
\end{equation}
(\ref{eq. 39}) states that different outer shifts of falling factorial transforms and also different inner shifts of inverse transform are orthogonal.

\subsection{Basis shifts}

In this section, our goal is to illustrate the treatment of falling factorial transform when this transform is applied to the terms $(x)_n (x)_m$ and $x^n x^m$. Hence, by (\ref{eq. 6}) we can write:
\begin{equation}
	\label{eq. 40}
	\begin{split}
		FFT^{-1} ((x)_n f(x)) & = e^{-x} \sum_{k=0}^{\infty} \frac{(k)_n f(k)}{k!} x^k = e^{-x} \sum_{k=0}^{\infty} \frac{f(k)}{(k-n)!} x^k \\
		& = e^{-x} \sum_{k=0}^{\infty} \frac{f(k+n)}{k!} x^{k+n} = x^n FFT^{-1} (f(x+n))
	\end{split}
\end{equation}
By applying $f(x) = FFT_t (g(t))(x-n)$ in (\ref{eq. 40}), we get:
\begin{equation}
	\label{eq. 41}
	FFT_t \Bigl(t^n g(t)\Bigr)(x) = (x)_n FFT_t \Bigl(g(t)\Bigr)(x-n)
\end{equation}
From (\ref{eq. 41}), (\ref{eq. 10}), and $a^x = \sum_{n \geq 0} (xlog(a))^n⁄n!$:
\begin{equation}
	\label{eq. 42}
	RFT(a^x f(x)) = a^{xE} RFT(f(x))
\end{equation}

\subsection{Binomial transform}

First of all, we define binomial transform for all positive real numbers as follows:
\begin{equation}
	\label{eq. 43}
	BT(f(x)) = \sum_{n=0}^{\infty} \binom{x}{n} f(n)
\end{equation}
Thus, we can use (\ref{eq. 1}) and (\ref{eq. 2}):
\begin{equation}
	\label{eq. 44}
	BT(f(x)) = FFT(\sum_{n=0}^{\infty} f(n) \frac{x^n}{n!}) = FFT(e^x FFT^{-1} (f(x)))
\end{equation}
For inverse binomial transform, we have:
\begin{equation}
	\label{eq. 45}
	f(x) = FFT(e^x FFT^{-1} (BT^{-1} (f(x))))
\end{equation}
\begin{equation}
	\label{eq. 46}
	BT^{-1} (f(x)) = FFT(e^{-x} FFT^{-1} (f(x)))
\end{equation}

\subsection{Convolution}

As we performed in the previous section, we define convolution (or binomial convolution \cite{ref. 3}) for positive real numbers as follows:
\begin{equation}
	\label{eq. 47}
	conv(f(x),g(x)) = \sum_{n=0}^{\infty} \binom{x}{n} f(x-n) g(n)
\end{equation}
By applying $FFT^{-1}$ on convolution:
\begin{equation}
	\label{eq. 48}
	\begin{split}
		& FFT^{-1} (conv(f(x),g(x))) = FFT^{-1} \Bigl(\sum_{k=0}^{\infty} \binom{x}{k} f(x-k) g(k) \Bigr) \\
		& = e^{-x} \sum_{n=0}^{\infty} \sum_{k=0}^{n} \binom{n}{k} f(n-k) g(k) \frac{x^k}{k!} = e^{-x} \Bigl(\sum_{n=0}^{\infty} \frac{f(n)}{n!} x^n \Bigr) \Bigl(\sum_{n=0}^{\infty} \frac{g(n)}{n!} x^n \Bigr) \\
		& = e^x FFT^{-1} (f(x)) FFT^{-1} (g(x))
	\end{split}
\end{equation}
In other expression:
\begin{equation}
	\label{eq. 49}
	conv(f(x),g(x)) = FFT(e^x FFT^{-1} (f(x)) FFT^{-1} (g(x)))
\end{equation}
We can consider $f(x) = 1$:
\begin{equation}
	\label{eq. 50}
	conv(1,g(x)) = FFT(e^x FFT^{-1} (g(x))) = BT(g(x))
\end{equation}
Which (\ref{eq. 50}) shows that binomial transform is a special case of convolution. Consecutive convolutions can easily be derived by (\ref{eq. 48}) as follows:
\begin{equation}
	\label{eq. 51}
	FFT^{-1} (conv^{(n)} (f(x))) = e^{nx} (FFT^{-1} (f(x)))^{n+1}
\end{equation}
\begin{equation}
	\label{eq. 52}
	f^n (x) = e^x FFT^{-1} (conv^{(n-1)} (FFT(e^{-x} f(x))))
\end{equation}
By (\ref{eq. 6}) and $f(x)=e^x FFT^{-1} (F(x))$:
\begin{equation}
	\label{eq. 53}
	\Bigl(\sum_{k=0}^{\infty} \frac{F(k) x^k}{k!}\Bigr)^{n+1} = \sum_{k=0}^{\infty} \frac{conv^{(n)} \Bigl(F(x)\Bigr)(k)}{k!} x^k
\end{equation}

\subsection{Scaling}

By utilizing $f(ax)=a^{x \partial_x} f(x)$ from \cite{ref. 2} we get:
\begin{equation}
	\label{eq. 54}
	\begin{split}
		FFT(f(ax)) & = FFT(a^{x \partial_x} f(x)) \\
		& = FFT \Bigl(\sum_{n=0}^{\infty} \sum_{k=0}^{n} {n \brace k}  \frac{(log(a))^n}{n!} x^k \partial_x^k f(x) \Bigr) \\
		& = \sum_{n=0}^{\infty} \sum_{k=0}^{n} {n \brace k} \frac{(log(a))^n}{n!} FFT(x^k \partial_x^k f(x)) \\
		& = \Bigl(\sum_{n=0}^{\infty} \sum_{k=0}^{n} {n \brace k} \frac{(log(a))^n}{n!} (x)_k \Bigl(\frac{\Delta_x}{E}\Bigr)^k \Bigr) FFT(f(x))
	\end{split}
\end{equation}
We know that $a^{x \Delta_x /E} = \sum_{n=0}^{\infty} \sum_{k=0}^{n} {n\brace k} \frac{(log(a))^n}{n!} (x)_k (\Delta_x/E)^k$, which can be proven by induction and the Taylor series. Thus, we obtain a formula for scaling property as follows:
\begin{equation}
	\label{eq. 55}
	FFT(f(ax)) = a^{x \frac{\Delta_x}{E}} FFT(f(x))
\end{equation}
Which $\Delta_x/E$ is equal to the backward difference.

\subsection{Multiplication}

From (\ref{eq. 38}) and (\ref{eq. 40}) we have:
\begin{equation}
	\label{eq. 56}
	FFT^{-1} ((x)_n (x)_m) = x^n FFT^{-1} ((x+n)_m) = x^n m! L_m^{(n-m)} (-x)
\end{equation}
Now using the definition of Laguerre polynomials:
\begin{equation}
	\label{eq. 57}
	(x)_n (x)_m = FFT(x^n m!L_m^{(n-m)} (-x)) = \sum_{k=0}^{m} \binom{n}{k} \binom{m}{k} k! (x)_{n+m-k}
\end{equation}
Hence:
\begin{equation}
	\label{eq. 58}
	\begin{split}
		& FFT^{-1} (f(x)g(x)) = FFT^{-1} \Bigl(\sum_{n \geq 0} \sum_{m \geq 0} \frac{\Delta^n f(0)}{n!}  \frac{\Delta^m g(0)}{m!} (x)_n (x)_m \Bigr) \\
		& = FFT^{-1} \Bigl(\sum_{n \geq 0} \sum_{m \geq 0} \frac{\Delta^n f(0)}{n!}  \frac{\Delta^m g(0)}{m!} \sum_{k=0}^{\infty} \binom{n}{k} \binom{m}{k} k! (x)_{n+m-k} \Bigr) \\
		& =\sum_{n \geq 0} \sum_{m \geq 0} \frac{\Delta^n f(0)}{n!}  \frac{\Delta^m g(0)}{m!} \sum_{k=0}^{\infty} \binom{n}{k} \binom{m}{k} k! x^{n+m-k} \\
		& = \sum_{n \geq 0} \sum_{m \geq 0} \sum_{k=0}^{\infty} \frac{\Delta^n f(0) x^{n-k}}{(n-k)!} \frac{\Delta^m g(0) x^{m-k}}{(m-k)!} \frac{x^k}{k!} \\
		& = \sum_{k=0}^{\infty} \partial_x^k (FFT^{-1} (f(x))) \partial_x^k (FFT^{-1} (g(x))) \frac{x^k}{k!}
	\end{split}
\end{equation}
The conflation of (\ref{eq. 58}) and (\ref{eq. 6}) gives the Hadamard product of the power series.

Let $f(x)=FFT(F(x))$ and $g(x)=FFT(G(x))$:
\begin{equation}
	\label{eq. 59}
	FFT(F(x))FFT(G(x)) = FFT(\sum_{k=0}^{\infty} \partial_x^k (f(x)) \partial_x^k (g(x)) \frac{x^k}{k!})
\end{equation}
By (\ref{eq. 46}) and (\ref{eq. 48}):
\begin{equation}
	\label{eq. 60}
	\begin{split}
		FFT^{-1} (BT^{-1} (conv(f(x),g(x)))) & = e^{-x} FFT^{-1} (conv(f(x),g(x))) \\
		& = FFT^{-1} (f(x)) FFT^{-1} (g(x))
	\end{split}
\end{equation}
Now let $f(x) = FFT(F(x))$ and $g(x) = FFT(G(x))$:
\begin{equation}
	\label{eq. 61}
	FFT(F(x)G(x))= BT^{-1} (conv(FFT(F(x)),FFT(G(x))))
\end{equation}

\subsection{Power series coefficients}

One of the useful applications of the falling factorial transform is finding the power series coefficient without getting involved in the derivative calculation of functions. The idea is simple and helps simplify problems. Applying inverse falling factorial transform on (\ref{eq. 44}), we get the following formulas:
\begin{equation}
	\label{eq. 62}
	FFT^{-1} (BT(\Gamma(x+1)a(x))) = e^x FFT^{-1} (\Gamma(x+1)a(x)) = \sum_{n=0}^{\infty} a_n x^n = f(x)
\end{equation}
\begin{equation}
	\label{eq. 63}
	a(x) = \frac{1}{\Gamma (x+1)} FFT(e^{-x} f(x))
\end{equation}

\subsection{Falling factorial transform with other transforms}

Here, we want to discuss the application of FFT on Borel transformation. In \cite{ref. 2}, Borel transformation is defined as:
\begin{equation}
	\label{eq. 64}
	B(f(x)) = \int_{0}^{\infty} e^{-t} f(xt) dt
\end{equation}
So, using $FFT$ and (\ref{eq. 55}):
\begin{equation}
	\label{eq. 65}
	\begin{split}
		FFT(B(f(x))) & = \int_{0}^{\infty} e^{-t} FFT(f(xt))dt = \int_{0}^{\infty} e^{-t} t^{x \frac{\Delta_x}{E}} FFT(f(x))dt \\
		& = \Gamma (x \frac{\Delta_x}{E} + 1) FFT(f(x))
	\end{split}
\end{equation}
Additionally, $FFT^{-1}$ can be applied to inverse Laplace transform. By the definition of both transforms \cite{ref. 4}:
\begin{equation}
	\label{eq. 66}
	\begin{split}
		& FFT^{-1} (\mathcal{L}^{-1} (f(x))) = FFT^{-1} (\frac{1}{2\pi i} \int_{\gamma - i\infty}^{\gamma +i\infty} f(s) e^{sx} ds) \\
		& = e^{-x} \sum_{n=0}^{\infty} \frac{x^n}{n!} \frac{1}{2\pi i}\int_{\gamma - i\infty}^{\gamma +i\infty} f(s) e^{sn} ds = e^{-x} \frac{1}{2\pi i}\int_{\gamma - i\infty}^{\gamma +i\infty} f(s) e^{xe^s} ds \\
		& = \frac{1}{2\pi i}\int_{\gamma - i\infty}^{\gamma +i\infty} f(s) e^{x(e^s-1)} ds
	\end{split}
\end{equation}
Also, we can apply RFT on Laplace transform \cite{ref. 6} using the Mellin transform properties:
\begin{equation}
	\label{eq. 67}
	\begin{split}
		RFT(e^x \mathcal{L}(f(x))) & = RFT \Bigl(\int_{0}^{\infty} f(t) e^{(1-t)x} dt \Bigr) = \int_{0}^{\infty} f(t)RFT(e^{(1-t)x})dt \\
		& = \int_{0}^{\infty} f(t) t^{-x} dt = \mathcal{M}(\frac{1}{x} f(\frac{1}{x})) = \mathcal{M}_t \Bigl(f(t) \Bigr)(-x-1)
	\end{split}
\end{equation}

\subsection{Fractional derivatives and differences}

Considering the definition of falling factorial transform discussed in section \ref{section 1}, we can calculate fractional derivatives of an arbitrary function as follows:
\begin{equation}
	\label{eq. 68}
	\partial_t^x (f(t) e^t) = e^t \sum_{n=0}^{\infty} \binom{x}{n} \partial_t^n f(t) = e^t FFT_x (f(x+t))
\end{equation}
So:
\begin{equation}
	\label{eq. 69}
	\partial_t^x (f(t)) = FFT_x (e^{-x} f(x+t))
\end{equation}
We can define fractional differences as follows:
\begin{equation}
	\label{eq. 70}
	\Delta_t^x (f(t)) = BT_x^{-1} (f(x+t)) = FFT_x (e^{-x} FFT_x^{-1} (f(x+t)))
\end{equation}
Some Interesting identities about fractional derivatives and differences can be obtained as below:
\begin{equation}
	\label{eq. 71}
	\begin{split}
		\partial_t^x (\partial_s^t (f(s) e^s)) & = e^s \partial_t^x FFT_t (f(t+s)) = e^s e^{s \Delta_t} \partial_t^x FFT(f(t)) \\
		& = e^s e^{s \Delta_t} FFT_x (e^{-x} FFT_t \Bigl(f(t)\Bigr)(x+t)) \\
		& = e^s e^{s \Delta_t} FFT_x (e^{-x} e^{t \partial_x} FFT_x (f(x))) \\
		& = e^s e^{s \Delta_t} e^t e^{t \Delta_x} FFT(e^{-x} FFT(f(x))) \\
		& = e^{s E_t} e^{t E_x} FFT(e^{-x} FFT(f(x)))
	\end{split}
\end{equation}
\begin{equation}
	\label{eq. 72}
	\begin{split}
		\Delta_t^x (\partial_s^t (f(s) e^s)) & = FFT_x (e^{-x} FFT_x^{-1} (e^s FFT_t \Bigl(f(t+s)\Bigr)(x+t))) \\
		& = FFT_x (e^{-x} FFT_x^{-1} (e^s (1+ \Delta_x)^t FFT_x (f(x+s)))) \\
		& = e^s FFT_x (e^{-x} (1+ \partial_x)^t f(x+s))
	\end{split}
\end{equation}
\begin{equation}
	\label{eq. 73}
	\begin{split}
		\partial_t^x (e^t \Delta_s^t (f(s))) & = e^t FFT_x \Bigl(FFT_t (e^{-t} FFT_t^{-1} (f(t+s))\Bigr)(x+t)) \\
		& = e^t FFT_x (e^{t \partial_x} FFT_x (e^{-x} FFT_x^{-1} (f(x+s)))) \\
		& = e^{t E_x} FFT_x^{(2)} (e^{-x} FFT_x^{-1} (f(x+s)))
	\end{split}
\end{equation}
\begin{equation}
	\label{eq. 74}
	\begin{split}
		& \Delta_t^x (\Delta_s^t (f(s))) \\
		& = FFT_x (e^{-x} FFT_x^{-1} (FFT_t \Bigl(e^{-t} FFT_t^{-1} (f(t+s)) \Bigr)(x+t))) \\
		& = FFT_x (e^{-x} (1+ \partial_x)^t FFT_x^{-1} (FFT_x (e^{-x} FFT_x^{-1} (f(x+s))))) \\
		& = FFT_x (e^{-x} (1+ \partial_x)^t (e^{-x} FFT_x^{-1} (f(x+s))))
	\end{split}
\end{equation}
To simplify the result, we prune the term $(1+ \partial_x )^t (e^{-x} FFT_x^{-1} (f(x+s)))$:
\begin{equation}
	\label{eq. 75}
	\begin{split}
		& (1+ \partial_x)^t (e^{-x} FFT_x^{-1} (f(x+s))) = BT_t (\partial_x^t (e^{-x} FFT_x^{-1} (f(x+s)))) \\
		& = e^{-x} BT_t (BT_t^{-1} (\partial_x^t (FFT_x^{-1} (f(x+s))))) \\
		& = e^{-x} \partial_x^t (FFT_x^{-1} (f(x+s)))
	\end{split}
\end{equation}
So, we get:
\begin{equation}
	\label{eq. 76}
	\Delta_t^x (\Delta_s^t (f(s))) = FFT_x (e^{-2x} \partial_x^t (FFT_x^{-1} (f(x+s))))
\end{equation}
Also, we can derive a formula for fractional integrals using the definition of RFT as follows:
\begin{equation}
	\label{eq. 77}
	\begin{split}
		\partial_t^{-x} (f(t)) & = \frac{1}{\Gamma (x)} \int_{t}^{\infty} f(u) u^{x-1} du = \frac{1}{\Gamma (x)} \int_{0}^{\infty} f(u+t) (u+t)^{x-1} du\\
		& = \frac{1}{\Gamma (x)} \sum_{n=0}^{\infty} \binom{x-1}{n} t^n \int_{0}^{\infty} f(u+t) u^{x-n-1} du\\
		& = \frac{1}{\Gamma (x)} \sum_{n=0}^{\infty} \frac{t^n}{n!} \mathcal{M}_x (\partial_x^n f(x+t)) = \frac{1}{\Gamma (x)} \mathcal{M}_x (f(x+2t)) \\
		& = RFT_x (f(x+2t)e^x)
	\end{split}
\end{equation}

\subsection{Other properties}

Another definition for inverse falling factorial transform can be derived as follows:
\begin{equation}
	\label{eq. 78}
	\begin{split}
		& FFT^{-1} (f(x)) = e^{-x} \sum_{n=0}^{\infty} \frac{E^n f(0) x^n}{n!} = e^{-x} \sum_{n=0}^{\infty} \frac{e^{n \partial_x} f(0) x^n}{n!} \\
		& = e^{-x} \sum_{n=0}^{\infty} \sum_{k=0}^{\infty} \frac{n^k \partial_x^k f(0) x^n}{n! k!} = e^{-x} \sum_{n=0}^{\infty} \sum_{k=0}^{\infty} \frac{\partial_x^k f(0) D_x^k x^n}{n! k!} = e^{-x} f(D_x) \{e^x\}
	\end{split}
\end{equation}
Applying FFT and (\ref{eq. 5}):
\begin{equation}
	\label{eq. 79}
	f(x) = FFT(e^{-x} f(D_x) \{e^x \}) = \partial_t^x (f(D_t) \{e^t \}) \Big\rvert_{t=0}
\end{equation}
It shows that there exists an elementary function (i.e., exponential function) and all other functions, indeed, are operators that operate on this elementary function.

Suppose $f(x) = (x)_{-n}$. By (\ref{eq. 6}) and the definition of incomplete gamma function \cite{ref. 5}:
\begin{equation}
	\label{eq. 80}
	\begin{split}
		FFT^{-1} ((x)_{-n}) & = e^{-x} \sum_{k=0}^{\infty} \frac{(k)_{-n} x^k}{k!} = e^{-x} \sum_{k=0}^{\infty} \frac{x^k}{(k+n)!} \\
		& = x^{-n} - \frac{\Gamma (n,x)}{(n-1)!} x^{-n}
	\end{split}
\end{equation}
One can prove by induction that $\Gamma(n,x) x^{-n} = (-1)^{n-1} \partial_x^{n-1} (e^{-x}/x)$. So, with (\ref{eq. 6}), $FFT(e^{1-x}/(1-x))=e \Gamma(x+1)$, (\ref{eq. 16}), and (\ref{eq. 33}) we have:
\begin{equation}
	\label{eq. 81}
	FFT^{-1} ((x)_{-n}) = x^{-n} + FFT^{-1} \Bigl(\frac{(-1)^{n-1}}{(n-1)!} \Delta^{n-1} e^{\Delta} (e \Gamma (x+1)) \Bigr)
\end{equation}
\begin{equation}
	\label{eq. 82}
	FFT(x^{-n}) = (x)_{-n} - \frac{(-1)^{n-1}}{(n-1)!} \Delta^{n-1} e^{\Delta} (e \Gamma (x+1))
\end{equation}
\begin{equation}
	\label{eq. 83}
	\begin{split}
		FFT(\mathcal{L} (f(x))) & = FFT(\sum_{n=1}^{\infty} f^{(n-1)} (0) x^{-n}) \\
		& = \sum_{n=1}^{\infty} f^{(n-1)} (0) (x)_{-n} - f(-\Delta) e^{\Delta} (e \Gamma (x+1))
	\end{split}
\end{equation}
Also for inverse RFT:
\begin{equation}
	\label{eq. 84}
	\begin{split}
		RFT^{-1} (x^{\overline{-n}}) & = e^x \sum_{k=0}^{\infty} \frac{(-1)^k (-k)^{\overline{-n}} x^k}{k!} = e^x \sum_{k=0}^{\infty} \frac{(-1)^{k+n} x^k}{(k+n)!} \\
		& = x^{-n} - \frac{\Gamma (n,-x)}{(n-1)!} x^{-n}
	\end{split}
\end{equation}
From (\ref{eq. 9}) we have $RFT(e^x) = \delta[x]$ (which $\delta [x] = \delta_{x,0}$). So, by RFT version of (\ref{eq. 41}):
\begin{equation}
	\label{eq. 85}
	\begin{split}
		& RFT(e^x) = x^{\overline{n}} RFT \Bigl(\frac{e^x}{x^n} \Bigr)(x+n) \\
		& \rightarrow RFT \Bigl(\frac{e^x}{x^n} \Bigr) = x^{\overline{-n}} \delta [x-n] = \frac{\delta [x-n]}{(x-1)_n}
	\end{split}
\end{equation}
We can rewrite this formula in the FFT form by (\ref{eq. 10}):
\begin{equation}
	\label{eq. 86}
	\begin{split}
		FFT \Bigl(\frac{e^{-x}}{x^n} \Bigr) = (-1)^n RFT \Bigl(\frac{e^x}{x^n} \Bigr)(-x) = (x)_{-n} \delta [x+n] = \frac{\delta [x+n]}{(x+1)^{\overline{n}}}
	\end{split}
\end{equation}

Let $G(x) = RFT^{-1} (g(x))$. So, by the multiplication property and (\ref{eq. 8}) we have:
\begin{equation}
	\label{eq. 87}
	\begin{split}
		g(x) \mathcal{M}(f(x)) & = g(x) \Gamma (x) RFT(e^x f(x)) = \Gamma (x) RFT(e^x f(x)) RFT(G(x)) \\
		& = \Gamma (x) RFT \Bigl(\sum_{k=0}^{\infty} \partial_x^k (e^x f(x)) \partial_x^k (G(x)) \frac{(-1)^k x^k}{k!} \Bigr) \\
		& = \mathcal{M} \Bigl(e^{-x} \sum_{k=0}^{\infty}\partial_x^k (e^x f(x)) \partial_x^k (G(x)) \frac{(-1)^k x^k}{k!} \Bigr) \\
		& = \mathcal{M} \Bigl(\sum_{k=0}^{\infty} \Bigl(\sum_{i=0}^{k} \binom{k}{i} f^{(i)} (x) \Bigr) \partial_x^k (G(x)) \frac{(-1)^k x^k}{k!} \Bigr) \\
		& = \mathcal{M} \Bigl(\sum_{i=0}^{\infty} \sum_{k=0}^{\infty} \frac{f^{(i)} (x) G^{(k+i)} (x)}{i!} \frac{(-x)^{k+i}}{k!} \Bigr) \\
		& = \mathcal{M} \Bigl( \sum_{i=0}^{\infty} \frac{f^{(i)} (x) G^{(i)} (x)}{i!} (-x)^i \Bigr)
	\end{split}
\end{equation}
We know that $g(x+n) = \sum_{k \geq 0} \frac{g^{(k)} (x)}{k!} n^k$, $T_k (-x) e^{-x} = \sum_{n \geq 0} \frac{(-1)^n n^k x^n}{n!}$. Hence, by basis shift property of the Mellin transform (from Table \ref{tab:Table 1}) and (\ref{eq. 8}):
\begin{equation}
	\label{eq. 88}
	\begin{split}
		\frac{1}{\Gamma (x)} e^{-E} g(x) \mathcal{M} (f(x)) & = \frac{1}{\Gamma (x)} \sum_{n \geq 0} \frac{(-1)^n}{n!} g(x+n) \mathcal{M}_t \Bigl(f(t) \Bigr) (x+n) \\
		& = \frac{1}{\Gamma (x)} \sum_{n \geq 0} \sum_{k \geq 0} \frac{(-1)^n n^k}{n!} \frac{g^{(k)} (x)}{k!} \mathcal{M}(x^n f(x)) \\
		& = \frac{1}{\Gamma (x)} \sum_{k \geq 0} \frac{g^{(k)} (x)}{k!} \mathcal{M} (T_k (-x) e^{-x} f(x)) \\
		& = \frac{1}{\Gamma (x)} \sum_{n \geq 0} \sum_{k \geq n} {k\brace n} \frac{g^{(k)} (x)}{k!} (-E)^n \mathcal{M} (e^{-x} f(x)) \\
		& = \frac{1}{\Gamma (x)} \sum_{n \geq 0} \frac{(e^{\partial_x}-1)^n g(x)}{n!} (-E)^n \mathcal{M} (e^{-x} f(x)) \\
		& = \sum_{n \geq 0} \frac{x^{\overline{n}} (-1)^n \Delta_x^n g(x)}{n!} RFT \Bigl(f(x) \Bigr)(x+n)
	\end{split}
\end{equation}
Using this formula, some properties of the Mellin transform from Table \ref{tab:Table 1} can be written in the RFT form. Table \ref{tab:Table 2}, shows these derived formulas.
\begin{table}[h!]
	\caption{Some properties of the RFT derived from the Mellin transform}
	\label{tab:Table 2}
	\centering
	\small
	\begin{tabular}{|c|c|}
		\hline
		\textbf{Function} & \textbf{RFT} \\
		\hline
		$x^a f(x)$ & $\frac{\Gamma (x+a)}{\Gamma (x)} RFT \Bigl(f(x) \Bigr)(x+a)$ \\
		\hline
		$f(ax)$ & $a^{-x} RFT \Bigl(e^{(1-\frac{1}{a})x} f(x) \Bigr)$ \\
		\hline
		$f(x^a)$ & $\frac{\Gamma (x/a)}{\lvert a \rvert \Gamma (x)} RFT \Bigl(e^{x-x^{\frac{1}{a}}} f(x) \Bigr) (\frac{x}{a})$ \\
		\hline
		$f(x)(log(x))^n$ & $\sum_{k=0}^{n} \binom{n}{k} \frac{\Gamma^{(k)} (x)}{\Gamma (x)} \partial_x^{n-k} (RFT(f(x)))$ \\
		\hline
		$(x\partial_x)^n f(x)$ & $\sum_{k=0}^{n} {n\brace k} x^{\overline{k}} \Delta_x^k RFT(f(x)) = (x\Delta_x)^n RFT(f(x))$ \\
		\hline
		$(\partial_x x)^n f(x)$ & $\sum_{k=0}^{n} {{n+1}\brace {k+1}} x^{\overline{k}} \Delta_x^k RFT(f(x))$ \\
		\hline
	\end{tabular}
\end{table}

The orthogonality property (i.e., equation (\ref{eq. 39})) and (\ref{eq. 38}) cause to emerge following identity:
\begin{equation}
	\label{eq. 89}
	\begin{split}
		& \sum_{k=0}^{\infty} \frac{(-x)^k}{k!} FFT^{-1} (f(x+k)) FFT^{-1} (g(x+k)) \\
		& = \sum_{n \geq 0} \sum_{m \geq 0} e^{-x} \Delta^n f(0) \frac{\Delta^m g(0)}{m!} (-x)^m \delta_{n,m} \\
		& = e^{-x} \sum_{k=0}^{\infty} \frac{(-1)^k \Delta^k f(0) \Delta^k g(0)}{k!} x^k
	\end{split}
\end{equation}

\section{Corollaries}
\label{section 3}

Table \ref{tab:Table 3} shows some functions that falling factorial transform applied to them. Also, the results of this table can be used for the RFT. A known series for the Riemann zeta function can be derived by the identities $\mathcal{M} ( \frac{1}{e^x-1}) = \Gamma (x) \zeta (x)$, $\frac{e^x}{e^x-1} = \sum_{n=-1}^{\infty} \frac{B_{n+1} (-1)^{n+1}}{(n+1)!} x^n$ (which $B_n$ are Bernoulli numbers), and (\ref{eq. 84}):
\begin{equation}
	\label{eq. 90}
	\begin{split}
		\zeta (x) & = RFT \Bigl(\frac{e^x}{e^x-1} \Bigr) = RFT \Bigl(\sum_{n=-1}^{\infty} \frac{B_{n+1} (-1)^{n+1}}{(n+1)!} x^n \Bigr) \\
		& = RFT \Bigl(\frac{e^x}{x} \Bigr) + \sum_{n=-1}^{\infty} \frac{B_{n+1} (-1)^{n+1}}{(n+1)!} x^{\overline{n}}
	\end{split}
\end{equation}
So, by (\ref{eq. 85}) and $B_0=1$:
\begin{equation}
	\label{eq. 91}
	\zeta (x) = \frac{\delta [x-1] - 1}{x-1} + \sum_{n=0}^{\infty} \frac{B_{n+1} (-1)^{n+1}}{(n+1)!} x^{\overline{n}}
\end{equation}

\begin{table}[h!]
    \caption{Functions and FFTs}
    \label{tab:Table 3}
    \centering
    \small
    \begin{tabular}{|c|c|p{50mm}|}
      \hline
      \textbf{Function} & \textbf{FFT} & \textbf{Proof} \\
      \hline
      $\frac{e^{-x}}{1-x}$ & $\Gamma (x+1)$ & Eq. (\ref{eq. 6}) \\
      \hline
      $\frac{\Gamma (y+1) e^{-x}}{(1-x)^{y+1}}$ & $\Gamma (x+y)$ & Eq. (\ref{eq. 6}) and $\Gamma (y+1)/(1-x)^{y+1} = \sum_{k \geq 0} (y+k)! x^k/k!$ \\
      \hline
      $x^s e^{-x}$ & $\Gamma (s+1) \delta [x-s]$ & $\partial_t^x (t^s) \Big\rvert_{t=0} = \Gamma (s+1) \delta[x-s]$ and Eq. (\ref{eq. 5}) \\
      \hline
      $x^n$ & $(x)_n$ & Definition of FFT \\
      \hline
      $(x)_n$ & $Z_n (x)$ & Section \ref{subsection 2.4} \\
      \hline
      $Z_n (x)$ & $\sum_{k=0}^{n} {n\brack k} (-1)^{n-k} Z_k (x)$ & Definition of $Z_n (x)$ \\
      \hline
      $T_n (x)$ & $x^n$ & Eq. (\ref{eq. 6}) and Dobinski formula \cite{ref. 1} \\
      \hline
      $\sum_{k=0}^{n} {n\brace k} T_k (x)$ & $T_n (x)$ & Definition of Touchard polynomials \\
      \hline
      $e^{(a-1)x}$ & $a^x$ & Definition of FFT and $a \geq 0$ \\
      \hline
      $e^{i \omega x}$ & $(\omega^{2}+1)^{x/2} e^{ix tan^{-1} (\omega)}$ & $a = 1+i\omega$ \\
      \hline
      $sin(\omega x)$ & $(\omega^{2}+1)^{x/2} sin(x tan^{-1} (\omega))$ & $sin(\omega x) = \frac{e^{i\omega x} - e^{-i\omega x}}{2i}$ \\
      \hline
      $sin(xtan(\omega))$ & $sin(\omega x)/cos^{x} (\omega)$ & $\omega \rightarrow tan(\omega)$ \\
      \hline
      $cos(\omega x)$ & $(\omega^{2}+1)^{x/2} cos(x tan^{-1} (\omega))$ & $cos(\omega x) = \frac{e^{i\omega x} + e^{-i\omega x}}{2}$ \\
      \hline
      $cos(xtan(\omega))$ & $cos(\omega x)/cos^{x} (\omega)$ & $\omega \rightarrow tan(\omega)$ \\
      \hline
      $x^{n+m} L_n^{(m)} (-x)$ & $(x)_{n+m} (x)_m /m!$ & Definition of FFT and Laguerre polynomials \\
      \hline
    \end{tabular}
\end{table}

\section{Conclusion}

In this paper, we discussed the properties of falling and rising factorial transforms and illustrated how FFT and RFT behave in different situations. Then, the FFT of some functions is calculated and shown in Table \ref{tab:Table 3}.

FFT and RFT are comprehensive tools that link power series into Newton series (i.e., derivatives into differences). These transforms have the following applications and properties: New classes of operators, orthogonality, new formulation for Hadamard product, finding power series coefficients of arbitrary functions without getting involved in the derivatives of functions, and fractional derivatives and differences.

\end{document}